\documentclass[11pt]{article}
\usepackage{amssymb}
\usepackage{amsmath}
\usepackage{tabularx}
\usepackage{graphicx}
\usepackage{epsfig,epsf,psfrag,epstopdf,subfigure}
\usepackage{color}
\usepackage{comment}
\usepackage{algorithm}
\usepackage{algorithmic}
\usepackage{url}
\usepackage{cleveref}

\newtheorem{lemma}{Lemma}

\newtheorem{definition}{Definition}
\newtheorem{remark}{Remark}

\newcommand{\be}{\begin{equation}}
\newcommand{\ee}{\end{equation}}
\newcommand{\ba}{\begin{aligned}}
\newcommand{\ea}{\end{aligned}}
\newcommand{\bea}{\begin{eqnarray}}
\newcommand{\eea}{\end{eqnarray}}

\def\qed{\hfill$\Box$}

\newcommand\tg{\tilde{\bf g}}

\newcommand\ci{\mathcal{I}}
\newcommand\cs{\mathcal{S}}

\newcommand\cv{\mathcal{V}}

\newcommand\vort{\boldsymbol{\omega}}
\newcommand\brho{\boldsymbol{\rho}}
\newcommand\bmu{\boldsymbol{\mu}}
\newcommand\bdelta{\boldsymbol{\delta}}
\newcommand\bm{\mathbf{m}}

\newcommand\bR{\mathbb{R}}

\newcommand\ec{e^{-c}}

\newcommand\bx{\mathbf{x}}
\newcommand\bJ{\mathbf{J}}

\newcommand\by{\mathbf{y}}

\newcommand\bxy{\bx-\by}

\newcommand\lbxy{\|\bxy\|}

\newcommand\dt{\Delta t}

\title{A new mixed potential representation for the 
equations of unsteady, incompressible flow} 

\author{Leslie Greengard%
  \thanks{Courant Institute of Mathematical Sciences,
    New York University, New York, New York 10012 and
    Flatiron Institute, Simons Foundation, New York, New York 10010.
    (Email: {greengard@courant.nyu.edu}).}
  \and
  Shidong Jiang%
  \thanks{Department of Mathematical Sciences,
    New Jersey Institute of Technology, Newark, New Jersey 07102
    This work was supported in part by the NSF under grant DMS-1720405
    and by the Flatiron Institute, a division of the Simons Foundation.
    (Email: {shidong.jiang@njit.edu}). 
  }
}

\begin{document}

\maketitle

\begin{abstract}
  We present a new integral representation for the unsteady, incompressible 
  Stokes or Navier-Stokes equations, based on
  a linear combination of heat and harmonic potentials. 
  For velocity boundary conditions, this leads to a coupled system of 
  integral equations: one for the normal component of velocity and one for the
  tangential components.
  Each individual equation is well-conditioned, and we show that using them
  in predictor-corrector fashion, combined with spectral deferred correction,
  leads to high-order accuracy solvers.
  The fundamental unknowns in the mixed potential representation are 
  densities supported on the boundary of the domain. We refer to one as the
  {\em vortex source}, the other as the {\em pressure source} and the
  coupled system as the {\em combined source integral equation}.  
\end{abstract}


\pagestyle{myheadings}
\thispagestyle{plain}
\markboth{L. Greengard and S. Jiang}
{Mixed potential representation for incompressible flow}

\section{Introduction} \label{intro}

We present a new integral representation for the numerical solution
of the unsteady, incompressible Navier-Stokes equations: 
\begin{equation}
  \frac{\partial {\bf u}}{\partial t}= \nu \Delta {\bf u}-\nabla p 
- ({\bf u} \cdot \nabla) {\bf u}  + {\bf f}, \qquad
   \nabla \cdot {\bf u}=0, 
  \label{nseq}
\end{equation}
or their linearization,
the unsteady Stokes equations,
\begin{equation}
  \frac{\partial {\bf u}}{\partial t}= \nu \Delta {\bf u}-\nabla p+{\bf F},
   \qquad \nabla \cdot {\bf u}=0.
  \label{stokeseq}
\end{equation}
The domain, which may be non-stationary, will be denoted at time 
$t \in [0,T]$ by $D(t)$ with boundary $\Gamma(t)$. The entire
space-time domain will be denoted by $D_T$ with boundary $\Gamma_T$.
Here, ${\bf u}(\bx,t)$ is the velocity field of interest, $\nu$ is the viscosity, and 
$p(\bx,t)$ is the pressure at a point $\bx \in D(t)$.
In \cref{nseq,stokeseq}, ${\bf f}$ and ${\bf F}$ are forcing
terms.
Initial conditions for the velocity are given by
\be
   {\bf u}(\bx,0) = {\bf u}_0(\bx), \qquad \bx \in D(0), \label{stokeseqs:ic}
\ee
and we restrict our attention to the case where 
``velocity" boundary conditions are prescribed:
\be
{\bf u}(\bx,t) = {\bf g}(\bx,t), \qquad (\bx,t) \in \Gamma_T.
\label{stokeseqs:dbc}
\ee

In this paper, we focus on the linearized problem, assuming that 
${\bf F}$ is known.  From a numerical perspective, it already contains the essential
difficulty faced by marching schemes for the Navier-Stokes equations,
which usually treat the nonlinear, advective term explicitly. 
That essential difficulty concerns computing the evolution of a 
diffusing velocity field, while maintaining the incompressibility condition
\begin{equation}
   \nabla \cdot {\bf u}=0
   \label{dfcond}
\end{equation}
through the addition of a pressure gradient.

Beginning with the work of Chorin and Temam \cite{chorin,temam}, 
one of the most popular approaches
for solving this problem is through the use of fractional step ``projection" methods. 
A simple version of such a scheme involves first solving
a diffusion equation for the velocity field with an explicit approximation of 
$\nabla p$ and ${\bf F}$, followed by the solution of a Poisson equation for the
pressure to enforce the incompressibility constraint.
Several decisions must be made in such schemes, 
including the choice of boundary conditions for the 
diffusion step and the choice of boundary conditions for the 
pressure correction/projection step.
We do not seek to review the literature here and refer the reader 
to \cite{BROWN,HENSHAW,liuliupego} for additional references and
a more thorough discussion.

To avoid fractional steps, an alternative is to use a gauge method.
Rather than solving the unsteady Stokes equations directly, one solves a system
of the form:
\begin{equation}
  \begin{aligned}
  \frac{\partial \bm}{\partial t} &= \nu \Delta \bm + {\bf F}, \\
   \Delta \phi &=  \nabla \cdot \bm,
  \end{aligned}
  \label{gauge}
\end{equation}
from which one obtains ${\bf u}$ and $p$ as 
\begin{eqnarray*}
   {\bf u} &=& \bm - \nabla \phi, \\
   p &=& \phi_t - \nu \Delta \phi.
\end{eqnarray*}
Such schemes require suitable boundary conditions for $\phi$ and $\bm$,
but avoid the fractional step and are more straightforward to discretize with 
high order accuracy in time (see, for example,
\cite{buttke,Cortez,eliu2003gauge,saye2017gauge,wang2000mcom}).

One can also obtain an unconstrained formulation by
taking the curl of the unsteady Stokes equations, yielding an equation for
the evolution of vorticity ${\vort} = \nabla \times {\bf u}$.
In three dimensions, we have
\begin{equation}
  \frac{\partial \vort}{\partial t}= \nu \Delta {\vort}+
\nabla \times {\bf F},
  \label{nseq2}
\end{equation}
while in two dimensions, 
\begin{equation}
  \frac{\partial \omega}{\partial t}= \nu \Delta {\omega}+
\left(\frac{\partial F_2}{\partial {x_1}} - \frac{\partial F_1}{\partial {x_2}} \right).
  \label{nseq2d}
\end{equation}
Here, ${\bf u} = (u_1,u_2)$, ${\bf F} = (F_1,F_2)$, and
vorticity is the scalar $\omega =
\frac{\partial u_2}{\partial {x_1}} - \frac{\partial u_1}{\partial {x_2}}$.
This approach is particularly natural in the two-dimensional setting, where
one can introduce a scalar stream function $\Psi$, with
\begin{equation} 
{\bf u} =  \nabla^{\perp} \Psi =  \left(\frac{\partial \Psi}{\partial x_2},
-\frac{\partial \Psi}{\partial x_1} \right), 
\label{gperp}
\end{equation}
so that the incompressibility constraint is automatically satisfied. 
It is easy to see that the stream function must satisfy the Poisson equation
\begin{equation} 
\Delta \Psi = -\omega.
\label{streampoisson}
\end{equation}
A major difficulty with this 
approach is that the boundary conditions for vorticity are nonlocal
\cite{anderson,Ben-Artzi2001,deangp,eliu2d,quartapelle}. 
Instead, one can also formulate the two-dimensional
Navier-Stokes equations entirely in terms of the stream function
\cite{Ben-Artzi2005,greengard1998sisc,Hou2009,jiang2013jcp}:
\begin{equation}
  \frac{\partial \Delta \Psi }{\partial t}  = \nu \Delta^2 {\Psi}-
\left(\frac{\partial F_2}{\partial x_1} - \frac{\partial F_1}{\partial x_2} \right).
  \label{nseq3}
\end{equation}
Since this is a fourth order partial differential equation,
one can directly impose velocity boundary conditions by 
specifying $\nabla^{\perp} \Psi$ 
on $\Gamma(t)$.
Unfortunately, the extension of this approach to three dimensions is 
much more complicated (see, for example, \cite{eliu3d}).

Finally, we should note that there is a Green's function for 
the linearized equations \eqref{stokeseq},
called the {\em unsteady Stokeslet}.
In \cite{jiang2012sisc},
integral equation methods were proposed using the corresponding 
layer potentials. While effective, they are somewhat complicated to 
implement with existing fast algorithms and high-order accurate quadrature
methods. We will return to this point in the concluding section. 

Here, we propose a new integral representation for the solution of the unsteady
Stokes equations that is divergence-free by construction, involves only the use
of harmonic and heat potentials, permits the natural imposition of velocity
boundary conditions, and is applicable in either two or three dimensions.
Since fast and high-order algorithms have been created for
harmonic and heat potentials over the past several decades, powerful numerical
machinery can immediately be brought to bear.
The heart of our approach is to find a particular solution to the 
inhomogeneous equation (accounting for the forcing term ${\bf F}$),
followed by a solution of the homogeneous, unsteady
Stokes problem to enforce the desired boundary conditions.
In three dimensions, the latter step involves a representation 
of the solution of the form
\[
\ba
{\bf u}(\bx,t)&=\nabla \phi(\bx,t)+ 
\nabla \times {\bf K}(\bx,t),\\
p(\bx,t)&=-\phi_t(\bx,t),
\ea
\]
where $\phi$ is harmonic and ${\bf K}$ satisfies the vector homogeneous 
heat equation. Both $\phi$ and ${\bf K}$ will be defined in terms of 
layer potentials on $\Gamma_T$, whose source densities will be referred 
to as the pressure source and vortex source, respectively.
Enforcing velocity boundary conditions will lead to 
the {\em combined source integral equation}.

The paper is organized as follows. 
In \cref{prelims},
we briefly summarize the necessary mathematical background.
In \cref{mixedpotsec}, we discuss the mixed potential representation
and derive the combined source integral equation.
In \cref{specsec},
we compute the spectrum and condition number of a fully implicit
version of the combined source integral
equations and in \cref{numsec}, we present numerical experiments. 
In \cref{altsec}, 
we investigate a kind of {predictor-corrector} scheme,
where we impose the normal and tangential boundary conditions sequentially.
In \cref{sdcsec}, we show how high-order 
accuracy can be achieved using a spectral deferred correction scheme. 
We conclude with an outline of future work.

\begin{remark}
For the sake of simplicity, we assume that $\nu = 1$ in the remainder of this paper.
This is easily accomplished
in the unsteady Stokes equations by rescaling the time variable.
\end{remark}

\section{Analytical Preliminaries} \label{prelims}

For a fixed Lipschitz domain $D$ in $\mathbb{R}^d$ with boundary $\Gamma$, we let
$L^2(D)$ denote the space of all square integrable functions in $D$ and we let
$L^2(\Gamma)$ denote the space of all square integrable functions on $\Gamma$. 
For the time-varying space-time cylinder $D_T \subset \mathbb{R}^d \times [0,T]$ 
with boundary $\Gamma_T$, we let
$L^2(D_T)$ denote the space of all square integrable functions in $D_T$ and we let
$L^2(\Gamma_T)$ denote the space of all square integrable functions on $\Gamma_T$. 
We briefly summarize the
necessary aspects of classical potential theory for
the Laplace and heat equations in $\mathbb{R}^d$ ($d=2,3$).

\subsection{Harmonic potentials}

The Green's function for the Laplace equation in free space is given by
\be\label{lgreen}
G_{\rm L}(\bx,\by)=\left\{\begin{array}{cc}
-\frac{1}{2\pi}\ln |\bx-\by| & \quad \text{in}\,\, \mathbb{R}^2,\\
\frac{1}{4\pi |\bx-\by|}  & \quad \text{in}\,\, \mathbb{R}^3.
\end{array}\right.
\ee

\begin{definition}
Let $\bx \in \bR^d$. 
The single layer potential $\mathcal{S}_{\rm L}$ with density $\rho \in L^2(\Gamma)$ 
is defined 
by
\be
\mathcal{S}_{\rm L}[\rho](\bx) = \int_{\Gamma}G_{\rm L}(\bx,\by)\rho(\by)ds(\by).
\ee 
The volume potential $\mathcal{V}_{\rm L}$ with density $f \in L^2(D)$ 
is defined 
by
\be
\mathcal{V}_{\rm L}[f](\bx) = \int_{D} G_{\rm L}(\bx,\by)f(\by)d\by.
\ee 
\end{definition}

\subsubsection{Jump relations}

For $\bx_0 \in \Gamma$,
the normal derivative of the single layer potential
$\mathcal{S}_{\rm L}[\rho]$ satisfies the jump relation
\be
\lim_{\epsilon\rightarrow 0+}
\frac{\partial \mathcal{S}_{\rm L}[\rho](\bx_0 \pm \epsilon \nu(\bx_0 ))}{\partial \nu(\bx_0)}
= \mp\frac{1}{2}\rho(\bx_0) + \mathcal{S}_{\rm L\nu}[\rho](\bx_0),
\ee 
where $\nu(\bx_0)$ is the unit outward normal vector to $\Gamma$ at the boundary point
$\bx_0$ and
\be
\mathcal{S}_{\rm L\nu}[\rho](\bx_0) = \mbox{p.v.}\int_{\Gamma}
\frac{\partial G_{\rm L}(\bx_0,\by)}{\partial \nu(\bx_0)} \rho(\by)ds(\by).
\ee
If $D$ is a Lipschitz domain, then
$\mathcal{S}_{\rm L}$ is a compact operator on
$L^2(\Gamma)$ and
$\mathcal{S}^\ast_{\rm L\nu}$ is bounded.
The latter is compact 
when $\Gamma$ is $C^1$ \cite{fabes,verchota}.

\subsubsection{Tangential derivatives}

In two dimensions, 
the tangential derivative of the single layer potential is denoted by
\be
\mathcal{S}_{\rm L\tau}[\rho](\bx_0)
= \int_{\Gamma} \frac{\partial G_{\rm L}(\bx_0,\by)}{\partial \tau(\bx_0)} \rho(\by)ds(\by),
\ee 
where $\tau(\bx_0)$ is the unit tangential vector at $\bx_0 \in \Gamma$.
$\mathcal{S}_{\rm L\tau}$ is defined in the Cauchy principal value sense and is
a bounded operator on $L^2(\Gamma)$. 
In three dimensions, the tangential derivatives can be written in the form
\be
\nu(\bx_0) \times \nabla \mathcal{S}_{\rm L}[\rho](\bx_0).
\ee 
This operator is, again, defined in the Cauchy principal value sense and 
bounded on $L^2(\Gamma)$. 
\subsection{Heat potentials}
The Green's function for the heat equation, $u_t = \Delta u$, is
\be
G_{\rm H}(\bx,t)
=\frac{1}{({4\pi t})^{d/2}}e^{-\frac{|\bx|^2}{4t}},
\quad \bx \in \mathbb{R}^d.
\ee

\begin{definition} \label{heatpotdefs}
Let $u_0 \in L^2(D(0))$, let
$f \in L^2(D_T)$ and let 
$\mu \in L^2(\Gamma_T)$.
Then, the initial heat potential $\ci_{\rm H}$ is defined by 
\begin{equation}\label{eq3.1}
\ci_{\rm H}[u_0](\bx,t) = \int_D G_{\rm H}(\bx-\by,t) u_0(\by)d\by,
\end{equation}
the volume heat potential $\cv_{\rm H}$ is defined by
\begin{equation}\label{eq3.2}
\cv_{\rm H}[g](\bx,t) = \int_0^t\int_{D(t')} G_{\rm H}(\bx-\by,t-t') f(\by,t')d\by dt',
\end{equation}
and the single layer heat potential $\cs_{\rm H}$ is defined by
\begin{equation}\label{eq3.3}
\cs_{\rm H}[\mu](\bx,t) = \int_0^t\int_{\Gamma(t')} G_{\rm H}(\bx-\by,t-t') \mu(\by,t')ds(\by)dt'.
\end{equation}
\end{definition}

\subsubsection{Jump relations}

It is well-known that
the initial heat potential $\ci_{\rm H}$ is a compact operator on
$L^2(D(0))$ and that the volume heat potential $\cv_{\rm H}$ is a
compact operator on $L^2(D_T)$.
As in the harmonic case, 
the normal derivative of the single layer heat potential 
$\cs_{\rm H\nu}[\mu]$
satisfies the jump relations
\begin{equation}
\lim_{\epsilon \rightarrow 0+}\cs_{\rm H\nu}[\mu](\bx_0\pm \epsilon \nu(\bx_0),t)
=\mp \frac{1}{2}\mu(\bx_0,t)+\cs_{\rm H\nu}[\mu](\bx_0,t), \quad \bx_0\in \Gamma,
\end{equation} 
where $\cs_{\rm H\nu}[\mu](\bx_0,t)$ is the principal value of 
\[
  \cs_{\rm H\nu}[\mu](\bx_0,t) = \int_0^t\int_{\Gamma(t')} \frac{\partial G_{\rm H}(\bx_0-\by,t-t')}
     {\partial \nu(\bx_0)} 
\mu(\by,t')ds(\by)dt'.
\]
In two dimensions, 
the tangential derivative of the single layer heat potential is denoted by
\be
\cs_{\rm H\tau}[\mu](\bx_0,t) = 
\int_0^t\int_{\Gamma(t')}
\frac{\partial G_{\rm H}(\bx_0-\by,t-t')}{\partial \tau(\bx_0)} \mu(\by,t')ds(\by)dt',
\ee 
where $\tau(\bx_0)$ is the unit tangential vector at $\bx_0 \in \Gamma$.
$\cs_{\rm H\tau}$ is defined in the Cauchy principal value sense and is
a bounded operator on $L^2(\Gamma)$. 

\subsubsection{The vector heat potential}

In three dimensions, we will also make use of a vector heat potential, defined by
\begin{equation}\label{vecheatpot}
{\bf K}_{\rm H}[\bJ](\bx,t) = \int_0^t\int_{\Gamma(t')} G_{\rm H}(\bx-\by,t-t') 
\bJ(\by,t')ds(\by)dt',
\end{equation}
where $\bJ$ is a tangential vector field on $\Gamma$.

\begin{definition}
For reasons that will become clear below,
we will refer to $\bJ$ or $\mu$ as the {\em vortex source}.
{\rm (}${\bf K}_{\rm H}[\bJ]$ will play a role analogous to that of the vector 
potential in electromagnetic theory, where the source is 
a surface electric current.{\rm )}
\end{definition}

It is straightforward to verify that, as in the electromagnetic case,
the tangential components of $\nabla \times {\bf K}_{\rm H}[\bJ]$ are given by
\cite{colton_kress_inverse}:
\begin{equation}
\lim_{\epsilon \rightarrow 0+}
\nu(\bx_0) \times \left( \nabla \times {\bf K}_{\rm H} \right) [\bJ]
(\bx_0\pm \epsilon \nu(\bx_0),t)
=  
\pm \frac{1}{2}\bJ(\bx_0,t) \, + \,
{\bf M}_{\rm H}[\bJ](\bx_0,t),
\label{mfiejump}
\end{equation} 
where
\begin{equation}\label{vecdble}
{\bf M}_{\rm H}[\bJ](\bx_0,t) = 
\nu(\bx_0) \times \left( \nabla \times {\bf K}_{\rm H} \right)[\bJ]
   (\bx_0,t), \quad \bx_0\in \Gamma,
\end{equation}
interpreted in the principal value sense. ${\bf M}_{\rm H}[\bJ]$ is a compact operator
when $\Gamma$ is $C^1$ \cite{colton_kress_inverse}.
Finally,
\begin{equation}\label{ndotK}
\nu(\bx_0) \cdot \nabla \times {\bf K}_{\rm H}[\bJ](\bx_0,t)
\end{equation} 
is defined in the Cauchy principal value sense, continuous across the boundary
of $\Gamma \in C^1$, and bounded on $L^2(\Gamma)$. 

\subsection{The Helmholtz decomposition of a vector field} \label{helmsec}

Let $D\subset \mathbb{R}^d$ be a bounded Lipschitz domain ($d = 2,3$).
It is well-known  \cite{girault}
that every vector field ${\bf F} \in L^2(D)$ has a decomposition of the form
\be
\label{hdecomp}
   {\bf F}=\nabla\phi+{\bf w},
\ee
where ${\bf w}$ is divergence-free (or {\em solenoidal})
and $\nabla\phi$ is curl-free (or {\em irrotational}).
We will sometimes write 
\be\label{hdecomp2}
   {\bf F}= {\bf F}_G + {\bf F}_S
\ee
instead of \eqref{hdecomp}, where ${\bf F}_G$ is irrotational and 
${\bf F}_S$ is solenoidal.

Without boundary conditions on  ${\bf w}$ or $\phi$, the Helmholtz decomposition is 
not unique. Nevertheless, assuming ${\bf F}$ is sufficiently smooth, a
simple explicit construction is easily computed.

\begin{lemma}\label{th:hdecomp} {\rm \cite{aris}}
  Let ${\bf F}$ be a twice differentiable vector field in a domain $D$
 with boundary $\Gamma$ in $\mathbb{R}^3$, and let
\[
 \phi(\bx) = -\int_D G_{\rm L}(\bx - \by) \left( \nabla_{\by} \cdot {\bf F}(\by) \right) \,d\by
+ \int_\Gamma G_{\rm L}(\bx - \by) \left( \nu(\by) \cdot {\bf F}(\by) \right) \, ds(\by),
\]
\[
 {\bf A}(\bx) = \int_D G_{\rm L}(\bx - \by)  \left( \nabla_{\by} \times {\bf F}(\by) \right) \,d\by
- \int_\Gamma G_{\rm L}(\bx - \by) \left( \nu(\by) \times {\bf F}(\by) \right) \, ds(\by).
\]
Then
\[ {\bf F} = \nabla \times {\bf A} + \nabla \phi. \]
In $\mathbb{R}^2$, if ${\bf F} = (F_1,F_2)$ is twice differentiable in a domain $D$
with boundary $\Gamma$, let $\phi$ be defined as above and let
\[
 \psi(\bx) = \int_D G_{\rm L}(\bx - \by)  \left( \nabla^\perp_{\by} \cdot {\bf F}(\by) \right) \,d\by
- \int_\Gamma G_{\rm L}(\bx - \by) \left( \tau(\by) \cdot {\bf F}(\by) \right) \, ds(\by),
\]
where $\nabla^\perp_{\bx} \cdot {\bf F}(\bx) = \frac{\partial F_2}{\partial x_1} - \frac{\partial F_1}{\partial x_2}$
and $\tau$ denotes the unit tangent vector along $\Gamma$.
Then
\[ {\bf F} = \nabla^\perp \psi + \nabla \phi. \]
\end{lemma}

Using the notation above, we can write this more compactly as
\begin{eqnarray}
\phi(\bx) &=& 
- \mathcal{V}_{\rm L}[\nabla \cdot {\bf F}](\bx) + \mathcal{S}_{\rm L}[\nu \cdot {\bf F}](\bx). 
\nonumber \\
\label{ppotrep}
{\bf A}(\bx) &=& 
\mathcal{V}_{\rm L}[\nabla \times {\bf F}](\bx) - \mathcal{S}_{\rm L}[\nu \times {\bf F}](\bx). \\
\psi(\bx) &=& 
\mathcal{V}_{\rm L}[\nabla^\perp \cdot {\bf F}](\bx) - \mathcal{S}_{\rm L}[\tau \cdot {\bf F}](\bx).
\nonumber
\end{eqnarray}
Both the harmonic volume potentials and the harmonic single layer potentials 
can be computed in optimal time, and 
with high order accuracy, using the fast multipole method 
and suitable quadrature rules
\cite{gimbutas2012jcp,fmm2,ethridge2001sisc,GLee,fmm,lprev,langston2011,biros2015cicp,
qbx3d}. 

\begin{remark}
In free space, there is an even simpler construction for the Helmholtz decomposition
(assuming sufficiently rapid decay of ${\bf F}$).

\begin{lemma}\label{th:hdecompfree} {\rm \cite{ladyzhenskaya}}
If ${\bf F} \in L^2(\bR^d)$, then
\[ {\bf F}_G = 
-\nabla \left( \nabla \cdot \int_{\bR^d} G_{\rm L}(\bx-\by) \, {\bf F}(\by) \, d\by \right), 
\qquad
{\bf F}_S = {\bf F} - {\bf F}_G,
\]
where $G_{\rm L}$ is the Green's function for the Laplace equation.
\end{lemma}
\end{remark}

\section{Potential theory for the unsteady Stokes equations} \label{mixedpotsec}

Before turning to the full boundary value problem, it is worth stating
a fundamental, but rarely used, fact about the unsteady Stokes equations
in the absence of physical boundaries.

\begin{lemma} \label{freespacelemma}
{\rm  (\cite{ladyzhenskaya}, chapter 4)}
Let ${\bf F}(\bx,t) \in L^2(\bR^d)$, where $d = 2,3$, with
the Helmholtz decomposition
\[ {\bf F}(\bx,t) = {\bf F}_S(\bx,t)  + {\bf F}_G(\bx,t), \]
where ${\bf F}_S$ is solenoidal and 
${\bf F}_G$ is irrotational. Then the solution to 
\eqref{stokeseq} 
in $\bR^d$ with divergence-free initial data ${\bf u}_0(\bx)$ is given by
\begin{equation} 
\begin{aligned}
{\bf u}^{(F)}(\bx,t) &= 
\ci_{\rm H}[{\bf u}_0](\bx,t) +
\cv_{\rm H}[{\bf F}_S](\bx,t) \\
\nabla p^{(F)}(\bx,t) &= {\bf F}_G(\bx,t),
\end{aligned}
\label{ufreesol}
\end{equation} 
where $G_{\rm H}(\bx,t)$ is the heat kernel.
(The operators $\ci_{\rm H}$ and $\cv_{\rm H}$ here are assumed to be defined 
on $\bR^d$ rather than a bounded domain $D$,)
\end{lemma}

In short, given the Helmholtz decomposition of the forcing term
${\bf F}$, the unsteady Stokes
equations have an explicit solution in free space by quadrature.
This turns out to be true in a bounded domain as well.

\begin{lemma}
Let ${\bf F}(\bx,t) \in L^2(D)$, where $d = 2,3$, with
the Helmholtz decomposition
\[ {\bf F}(\bx,t) = {\bf F}_S(\bx,t)  + \nabla \phi(\bx,t), \]
where ${\bf F}_S = \nabla \times {\bf A}$ in $\bR^3$ and
${\bf F}_S = \nabla^\perp \psi$ in $\bR^2$.
Then, a {\underline {\em particular}} solution to 
\cref{stokeseq,stokeseqs:ic} is given by
\begin{equation} 
\begin{aligned}
{\bf u}^{(F)}(\bx,t) &= 
\ci_{\rm H}[{\bf u}_0](\bx,t) +
\nabla \times \cv_{\rm H}[{\bf A}](\bx,t) \qquad {\rm in\ } \bR^3 \\
\nabla p^{(F)}(\bx,t) &= \nabla \phi(\bx,t), \\
& \\
{\bf u}^{(F)}(\bx,t) &= 
\ci_{\rm H}[{\bf u}_0](\bx,t) +
\nabla^\perp \cv_{\rm H}[\psi](\bx,t) \qquad {\rm in\ } \bR^2 \\
\nabla p^{(F)}(\bx,t) &= \nabla \phi(\bx,t),
\end{aligned}
\label{uparticular}
\end{equation} 
where the initial and volume heat potentials are given in
\cref{heatpotdefs}.
\end{lemma}

\noindent {\em Proof:}
It is straightforward to verify that the 
partial differential equation \eqref{stokeseq} is satisfied.
The fact that ${\bf u}^{(F)}(\bx,t)$ is divergence-free follows
immediately from \cref{freespacelemma} 
for the term $\ci_{\rm H}[{\bf u}_0](\bx,t)$ and by construction for the term
involving $\cv_{\rm H}[{\bf A}](\bx,t)$ or 
$\cv_{\rm H}[\psi](\bx,t)$.
\qed

Thus, from the preceding Lemma, we may represent the solution to 
the full unsteady Stokes equations in the form
\[ {\bf u} = {\bf u}^{(F)} + {\bf u}^{(B)}, \quad
\nabla p = \nabla p^{(F)} + \nabla p^{(B)}, \]
where
$({\bf u}^{(B)},\nabla p^{(B)})$ satisfy the homogeneous, linearized equations
\be\label{eq4.2}
\ba
  \frac{\partial {\bf u}^{(B)}}{\partial t}&=\Delta {\bf u}^{(B)}-\nabla p^{(B)},
 \qquad (\bx,t) \in D_T,  \\
   \nabla \cdot {\bf u}^{(B)}&=0,   \qquad (\bx,t) \in D_T,\\
   {\bf u}^{(B)}(\bx,0) &= {\bf 0}, \qquad \bx \in D(0), \\
   {\bf u}^{(B)}(\bx,t)
   &=\tg(\bx,t):= {\bf g}(\bx,t)-{\bf u}^{(F)}(\bx,t),
   \qquad (\bx,t) \in \Gamma_T.
\ea
\ee

There is a significant advantage in solving 
the homogeneous equations \eqref{eq4.2} rather than 
\cref{stokeseq,stokeseqs:ic,stokeseqs:dbc},
as we shall now see.

\subsection{The mixed potential representation}

Let us represent the solution to the homogeneous system,
\[ ({\bf u}^{(B)}(\bx,t),p^{(B)}(\bx,t)), \]
in terms of harmonic and heat layer potentials. In three dimensions, we define

\be\label{ierep3d}
\ba
{\bf u}^{(B)}(\bx,t)&=\nabla \cs_{\rm L}[\rho](\bx,t)+ 
\nabla \times {\bf K}_{\rm H}[\bJ](\bx,t),\\
p^{(B)}(\bx,t)&=- \frac{\partial}{\partial t} \cs_{\rm L}[\rho](\bx,t),
\ea
\ee
while in two dimensions, we define
\be\label{ierep}
\ba
{\bf u}^{(B)}(\bx,t)&=\nabla \cs_{\rm L}[\rho](\bx,t)+ 
\nabla^{\perp} \cs_{\rm H}[\mu](\bx,t),\\
p^{(B)}(\bx,t)&=- \frac{\partial}{\partial t} \cs_{\rm L}[\rho](\bx,t).
\ea
\ee
Here, $\rho$, $\bJ$ and $\mu$ are unknown 
{\em boundary densities} to be determined.
It is straightforward to verify that the representations
\eqref{ierep3d} and \eqref{ierep}
satisfy the first three equations in \eqref{eq4.2}.

\begin{definition}
Because of the preceding relations, we will refer to $\rho$ as 
the {\em pressure source} or {\em pressure source density}.
\end{definition}

\subsubsection{The combined source integral equation}

If we decompose the velocity field into a sum
of normal and tangential components on the boundary, then imposing velocity 
boundary conditions
leads, in two dimensions, to the system of integral equations
\be\label{bie}
\ba
\frac{1}{2}\rho(\bx,t)+\cs_{\rm L\nu}[\rho](\bx,t)+\cs_{\rm H\tau}[\mu](\bx,t)&=
\nu \cdot \tg(\bx,t),\\
\frac{1}{2}\mu(\bx,t)+\cs_{\rm H\nu}[\mu](\bx,t) - \cs_{\rm L\tau}[\rho](\bx,t) &=
-\tau \cdot \tg(\bx,t)
\ea
\ee
for the unknowns $\rho$ and $\mu$, where $\bx$ 
is a point on the boundary $\Gamma(t)$. 

In three dimensions, we obtain system of integral equations
\be\label{bie3d}
\ba
\frac{1}{2}\rho(\bx,t)+\cs_{\rm L\nu}[\rho](\bx,t)+
\nu(\bx,t) \cdot \nabla \times {\bf K}_{\rm H}[\bJ](\bx,t)
&= \nu(\bx,t) \cdot \tg(\bx,t),\\
\frac{1}{2}\bJ(\bx,t)+
{\bf M}_{\rm H}[\bJ](\bx,t)+
\nu(\bx,t) \times \nabla \cs_{\rm L}[\rho](\bx,t)
&= \nu(\bx,t) \times \tg(\bx,t).
\ea
\ee
for the unknowns $\rho$ and $\bJ$, where $\bx$ 
is a point on the boundary $\Gamma(t)$. 
We will refer to either \eqref{bie} or \eqref{bie3d} as the 
{\em combined source integral equation}.

One major advantage of the mixed potential representation is that
the unknown densities 
(the pressure source and the vortex source)
correspond to physical quantities of interest. 
The harmonic potential $\phi(\bx,t) = \cs_{\rm L}[\rho](\bx,t)$ 
determines the pressure, according to \cref{ierep3d,ierep},
while the heat potentials determine the vorticity.
More precisely, in three dimensions, 
\[
\vort(\bx,t) = \nabla \times \left(\nabla \times {\bf K}_{\rm H}\right)[\bJ](\bx,t) 
+ \nabla \times {\bf u}^{(F)},
\]
while in two dimensions,
\[
\omega(\bx,t) = -\partial_t \cs_{\rm H}[\mu](\bx,t) + 
\frac{\partial u_2^{(F)}}{\partial x_1} - \frac{\partial u_1^{(F)}}{\partial x_2}, 
\]
where ${\bf u}^{(F)} = (u_1^{(F)}, u_2^{(F)})$).
These relations may be of some direct interest in analysis.

\begin{remark}
It is, perhaps, worth
noting that in the mixed potential representation,
the boundary conditions for ${\bf K}$ and $\phi$ are local, but they yield
exact, nonlocal expressions for the pressure and vorticity through the formulae above.
\end{remark}

\subsection{Discretization} \label{discret-sec}

For the sake of simplicity, we restrict our attention to 
the integral equation system \eqref{bie} in two dimensions, and begin by 
semi-discretization in time (i.e., discretization 
with respect to the time variable alone). For this, we let 
\[ \rho_j=\rho(\bx,j\Delta t), \quad  \mu_j=\mu(\bx,j\Delta t). \]
\[ \brho_j = [\rho_0,\dots,\rho_j], \quad \bmu_j = [\mu_0,\dots,\mu_j] \]
We then write 
\[ \cs_{\rm H\nu}[\bmu_j](\bx,t) =  
\cs_{\rm H\nu}^{\rm far}[\bmu_{j-1}](\bx,t) +
\cs_{\rm H\nu}^{\rm loc}[\bmu_j](\bx,t), 
\]
\[ \cs_{\rm H\tau}[\bmu_j](\bx,t) =  
\cs_{\rm H\tau}^{\rm far}[\bmu_{j-1}](\bx,t) +
\cs_{\rm H\tau}^{\rm loc}[\bmu_j](\bx,t), 
\]
to denote the semi-discrete approximations of $\cs_{\rm H\nu}[\mu]$ and $\cs_{\rm H\tau}[\mu]$,
where
\be
\ba
&\cs_{\rm H\nu}^{\rm far}[\bmu_{j-1}](\bx,t) = \\ 
&\hspace{.4in} \sum_{l = 1}^{j-1} 
\int_{(l-1)\Delta t}^{l\Delta t} \int_{\Gamma(t')}
\frac{\partial G_{\rm H}(\bx-\by,t-t')}{\partial \nu(\bx)} 
P^{I}_{k}[\bmu_l](\by,t') ds(\by)dt'  ]\, , \\
&\cs_{\rm H\tau}^{\rm far}[\bmu_{j-1}](\bx,t) = \\
&\hspace{.4in} \sum_{l = 1}^{j-1} 
\int_{(l-1)\Delta t}^{l\Delta t} \int_{\Gamma(t')}
\frac{\partial G_{\rm H}(\bx-\by,t-t')}{\partial \tau(\bx)} 
P^{I}_{k}[ \bmu_l ](\by,t') ds(\by)dt' \, , \\
&\cs_{\rm H\nu}^{\rm loc}[\bmu_j](\bx,t) = \\  
&\hspace{.4in} \int_{(j-1)\Delta t}^{j\Delta t} \int_{\Gamma(t')}
\frac{\partial G_{\rm H}(\bx-\by,t-t')}{\partial \nu(\bx)} 
P^{I}_{k}[\bmu_j](\by,t') ds(\by)dt' \, ,  \\
&\cs_{\rm H\tau}^{\rm loc}[\bmu_j](\bx,t) = \\
&\hspace{.4in} \int_{(j-1)\Delta t}^{j\Delta t} \int_{\Gamma(t')}
\frac{\partial G_{\rm H}(\bx-\by,t-t')}{\partial \tau(\bx)} 
P^{I}_{k}[ \bmu_j ](\by,t') ds(\by)dt' \, .
\ea 
\ee 
Here, $P^{I}_{k}[\bmu_l]$
is the $k$th order Lagrange interpolant of the data 
\[ \{ \mu_l,\mu_{l-1},\dots,\mu_{l-k} \} \]
at the $(k+1)$ uniformly spaced
time points $\{ l \Delta t, (l-1) \Delta t, \dots, (l-k) \Delta t  \}$.

Note that we have separated out the contributions to 
$\cs_{\rm H\nu}$ and $\cs_{\rm H\tau}$ from 
the early time steps
($\cs_{\rm H\nu}^{\rm far}, \cs_{\rm H\tau}^{\rm far}$) from the contributions on the
most recent time interval
($\cs_{\rm H\nu}^{\rm loc}, \cs_{\rm H\tau}^{\rm loc}$). 
The superscript $I$ in the expression
$P^{I}_{k}[\bmu_l]$ indicates that
the latest time point $t_l = l \Delta t$ is being used in the polynomial interpolant.
Since $\cs_{\rm H\nu}^{\rm loc}$ and $\cs_{\rm H\tau}^{\rm loc}$ are linear operators acting on the 
densities $\{ \mu_j, \mu_{j-1},\dots,\mu_{j-k} \}$, we will also have occasion to write
\be
\ba
&\cs_{\rm H\nu}^{\rm loc}[\bmu_j](\bx, j \Delta t) =  \\
&\hspace{.4in} \left[ A_{\rm H\nu}^{j,k} \mu_{j-k} + \dots
 A_{\rm H\nu}^{j,1} \mu_{j-1} \right] + A_{\rm H\nu}^{j,0} \mu_{j} 
= F_{\rm H\nu}^{j,k}[\bmu_{j-1}] + A_{\rm H\nu}^{j,0} \mu_{j} 
 \\
&\cs_{\rm H\tau}^{\rm loc}[\bmu_j](\bx, j ]\Delta t) = \\
&\hspace{.4in} \left[ A_{\rm H\tau}^{j,k} \mu_{j-k} + \dots
 A_{\rm H\tau}^{j,1} \mu_{j-1} \right] + A_{\rm H\tau}^{j,0} \mu_{j}
= F_{\rm H\tau}^{j,k}[\bmu_{j-1}] + A_{\rm H\tau}^{j,0} \mu_{j}. 
\ea 
\ee 
This makes explicit the contributions of the densities at the various
time steps to the local potentials $\cs_{\rm H\nu}^{\rm loc}$ and $\cs_{\rm H\tau}^{\rm loc}$.

In the present paper, following previous work with the unsteady Stokeslet
\cite{jiang2012sisc}, we interchange the order of integration in space and time,
as proposed in \cite{li2009sisc}, and carry out the time integration analytically.
For $k = 1$, using an implicit interpolation rule, $P^{I}_1$, to achieve second order
accuracy in $\Delta t$,
the kernels of the spatial 
operators $A_{\rm H\nu}^{j,0}$  and $A_{\rm H\tau}^{j,0}$ are 
\be  \label{kernels2dheat}
\ba
G^{\rm loc}_{\rm H\nu}(\bx,\by)&=-\frac{(\bxy)\cdot\nu}{2\pi \lbxy^2}e^{-\frac{\lbxy^2}{4\Delta t}}+\frac{(\bxy)\cdot\nu}
{8\pi\Delta t}E_1\left(\frac{\lbxy^2}{4\Delta t}\right),\\
G^{\rm loc}_{\rm H\tau}(\bx,\by)&=-\frac{(\bxy)\cdot\tau}{2\pi \lbxy^2}e^{-\frac{\lbxy^2}{4\Delta t}}+\frac{(\bxy)\cdot\tau}
{8\pi\Delta t}E_1\left(\frac{\lbxy^2}{4\Delta t}\right),
\ea
\ee
where $E_1(x)=\int_x^\infty \frac{e^{-t}}{t} dt$ is the exponential integral function
\cite{nisthandbook}.

We will have occasion to make use of the explicit form of the interpolant as well.

\begin{definition}
$P^{E}_{k}[\bmu_{j-1}](t)$ is defined to be the 
$k$th order Lagrange {\em extrapolant} 
of the data 
\[ [\mu_{j-1},\dots,\mu_{j-(k+1)}]  \]
evaluate at $t = j \Delta t$.
\end{definition}

As in multistep methods for ordinary differential equations, when
the interpolation order $k > 1$, some
care is required in initialization - that is, computing the first $k-1$ time steps
with sufficient accuracy.
We will ignore this issue for the moment to avoid distractions.

Finally, the spatial integrals in 
$\cs_{\rm H\tau}^{\rm loc}$, $\cs_{\rm H\nu}^{\rm loc}$, 
$\cs_{\rm L\tau}$, and $\cs_{\rm L\nu}$ involve either logarithmic singularities
or principal value-type integrals. We use the quadrature schemes of 
\cite{alpert1999sisc} to discretize these integrals to sixteenth order accuracy
on smooth curves.

\subsection{History dependence and fast algorithms}

From a practical perspective, 
$\cs_{\rm H\tau}^{\rm far}$ and $\cs_{\rm H\nu}^{\rm far}$ clearly depend on the entire 
space-time history of the problem at hand. In the absence of suitable
algorithms, the cost of their evaluation would be prohibitive.
Fortunately, a number of fast algorithms have been developed for precisely this
purpose \cite{greengard2000acha,greengard1990cpam,lubichheat1,tauschheat1}
that permit their evaluation in $O(NM \log M)$ work rather than $O(N^2M^2)$ work, where
$N$ denotes the number of time steps and $M$ denotes the number of points in the 
discretization of the boundary.
We refer the reader to those papers for further details.
\section{Spectrum of the fully implicit combined source integral equation} 
\label{specsec}

For the sake of simplicity, we restrict our attention to 
the coupled integral equations \eqref{bie} in two dimensions on a stationary boundary.
While at first glance, this might appear to involve a compact perturbation of the 
identity, that is not the case.
Informally, this can be seen as follows. First, we note that
the operators $S_{\rm L\tau}$ and $S_{\rm H\tau}$ 
are compact perturbations of the operator $-\frac{1}{2}H$,
where $H$ denotes the Hilbert transform operator on the circle 
with perimeter $L$. $L$ here is the length of $\Gamma$ \cite{kress2014}:
\[ H[f](s) = \frac{1}{2 \pi} p.v. \int_0^{L} \cot \left( \frac{\pi (s-s')}{L} \right) 
f(s') ds'.
\]
Thus, the system can be written in the form 
\be
\begin{bmatrix}\frac{1}{2}I & \frac{1}{2}H\\
  -\frac{1}{2}H & \frac{1}{2}I
\end{bmatrix} \begin{bmatrix} \rho \\ \mu \end{bmatrix} + 
{\bf C}
\begin{bmatrix} \rho \\ \mu \end{bmatrix}  = 
\begin{bmatrix} 
\nu \cdot \tg(\bx,t)\\
-\tau \cdot \tg(\bx,t)
\end{bmatrix},
\ee
where ${\bf C}$ is compact.
Since $H^2=-I$, the determinant of the leading part of the system vanishes,
so that the coupled system fails to be a Fredholm equation of the second kind.

We now study the spectrum of the integral equation in detail when $\Gamma$ 
is a circle of radius $r$.
The resulting properties follow qualitatively for any smooth curve.
More specifically,
let us assume we are using a second-order accurate one-step implicit marching scheme,
as described in \cref{discret-sec}.
This yields an integral equation at the $j$th time step of the form
\be \label{syscircle}
\begin{bmatrix}\frac{1}{2}I+ \cs_{\rm L\nu} & A^{j,0}_{\rm H\tau} \\
  -\cs_{\rm L\tau}& \frac{1}{2}I+A^{j,0}_{\rm H\nu}\end{bmatrix}
\begin{bmatrix}\rho_j\\ \mu_j\end{bmatrix}
  =\begin{bmatrix}
\nu \cdot \tg - \cs^{\rm far}_{\rm H\tau}[\bmu_{j-1}] - F_{\rm H\tau}^{j,1}[\bmu_{j-1}] \\
-\tau \cdot \tg- \cs^{\rm far}_{\rm H\nu}[\bmu_{j-1}]] - F_{\rm H\nu}^{j,1}[\bmu_{j-1}]
\end{bmatrix}
\ee
using the notation of \cref{discret-sec}.

The kernels of the operators on the left-hand side of \eqref{syscircle} are given by 
\be
G_{\rm L\nu}(\bx,\by)=-\frac{(\bxy)\cdot\nu}{2\pi\lbxy^2},
\quad G_{\rm L\tau}(\bx,\by)=-\frac{(\bxy)\cdot\tau}{2\pi\lbxy^2}
\ee
and \eqref{kernels2dheat}.
On a circle of radius $r$, we have 
\[
\ba
\bx =(r\cos s',r\sin s'), &\quad \by=(r\cos s,r\sin s) \\
(\bxy)\cdot\nu &= r(1-\cos(s'-s)), &\quad (\bxy)\cdot\tau = r\sin(s'-s), \\
\lbxy^2 &= 2r^2(1-\cos(s'-s)).
\ea
\] 
We then have
\be
G_{\rm L\nu}(\bx,\by)=-\frac{1}{4\pi r},\quad
G_{\rm L\tau}(\bx,\by)=-\frac{1}{4\pi r}\cot\left(\frac{s'-s}{2}\right).
\ee
That is, the kernel of $S_{\rm L\nu}$ is constant and 
$S_{\rm L\tau}=-\frac{1}{2}H$ where $H$ is
the Hilbert transform on the unit circle. 
It is easy to verify that
all of these operators are diagonalized by the Fourier transform.
Thus, we only need to consider
the $2\times 2$ block for each Fourier mode $e^{iks}$ with $k \in  \mathbb{Z}$.

For $k=0$, we have
\be
\ba
\left(\frac{1}{2}I+S_{\rm L\nu}\right)[1] 
&= \frac{1}{2}-\frac{1}{4\pi r}\int_0^{2\pi}1\cdot rds=0,\\
S_{\rm L\tau}[1]&=-\frac{1}{2}H[1]=0,\\
\ea
\ee

$A^{j,0}_{\rm H\tau}[1]=0$ by symmetry, and
\be
\left(\frac{1}{2}I+A^{j,0}_{\rm H\nu}\right)[1]=\lambda_0
\ee
where 
\be
\lambda_0 =\frac{1}{2}-
\int_0^{2\pi} \left[  \frac{e^{-r^2(1-\cos s)/(2\Delta t)}}{4\pi} -
\frac{(1-\cos s) \, r^2 \, E_1\left(\frac{r^2(1-\cos s)}{2\Delta t}\right)}
{8\pi\Delta t} \right] ds.
\ee
Since $E_1(x)>0$ for $x>0$, we have $\lambda_0>0$ for any $r>0$ and $\Delta t>0$.
Thus, the system \eqref{syscircle} has 
eigenvalues $0$ and $\lambda_0$ with eigenvectors
$[1\,\,\, 0]^T$ and
$[0\,\,\, 1]^T$.

\begin{figure}[ht]
\centering
\includegraphics[width=.50\textwidth]{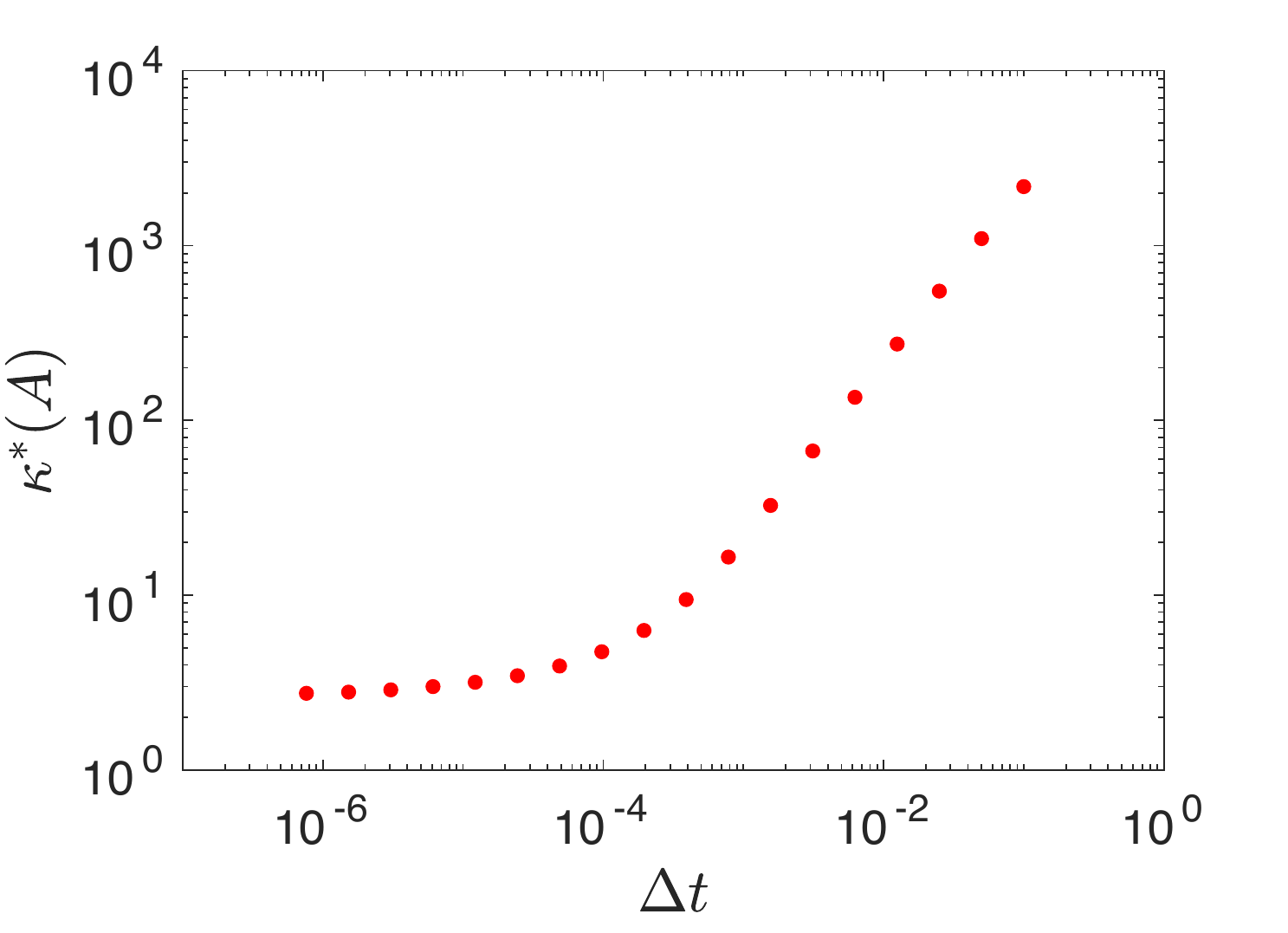}
\caption{The condition number of the fully implicit combined source
integral equation 
(in the complement of the one-dimensional nullspace) as a function of the 
time step for a circle of radius $r=0.6$, discretized with $128$ points.}
\label{condfig}
\end{figure}

For $k\neq 0$, we have
\be
\ba
\left(\frac{1}{2}I+S_{\rm L\nu}\right)[e^{iks}](s') &= \frac{1}{2}e^{iks'},\\
S_{\rm L\tau}[e^{iks}](s') &=-\frac{1}{2}H[e^{iks}](s')=\frac{1}{2}i\operatorname{sgn}(k) e^{iks'},\\
\left(\frac{1}{2}I+A^{j,0}_{\rm H\nu}\right)[e^{iks}](s')&=a_ke^{iks'},\\
A^{j,0}_{\rm H\tau}[e^{iks}](s') &= b_ke^{iks'},
\ea
\ee
with $a_k$, $b_k$ defined by the formulas
\be
\ba
a_k &= \frac{1}{2}-\frac{1}{4\pi}\int_{-\pi}^{\pi}e^{-\frac{r^2\sin^2(s/2)}{\Delta t}}\cos(ks)ds\\
&+\frac{r^2}{8\pi\Delta t}\int_{-\pi}^{\pi}(1-\cos s)E_1\left(\frac{r^2\sin^2(s/2)}{\Delta t}\right)\cos(ks)ds,\\
b_k &= \frac{i}{4\pi}\int_{-\pi}^{\pi}\cot\left(\frac{s}{2}\right)e^{-\frac{r^2\sin^2(s/2)}{\Delta t}}\sin(ks)ds\\
&-\frac{ir^2}{8\pi\Delta t}\int_{-\pi}^{\pi}\sin(s)E_1\left(\frac{r^2\sin^2(s/2)}{\Delta t}\right)\sin(ks)ds.
\ea
\ee
We note that $E_1(x)$ has a series expansion $E_1(x)=\ln x+\gamma+x+\frac{x^2}{4}+\ldots$
and that the spectrum of an integral operator with a smooth kernel decays exponentially 
fast. Thus, we have
\be
\ba
a_k &\approx \frac{1}{2}
+\frac{r^2}{8\pi\Delta t}\int_{-\pi}^{\pi}(1-\cos s)\ln\left(\sin^2(s/2)\right)\cos(ks)ds,\\
b_k &\approx \frac{i}{4\pi}\int_{-\pi}^{\pi}\cot\left(\frac{s}{2}\right)\sin(ks)ds \\
&\hspace{.5in} -\frac{ir^2}{8\pi\Delta t}\int_{-\pi}^{\pi}\sin(s)\ln\left(\sin^2(s/2)\right)\sin(ks)ds.
\ea
\ee
Using the facts (see, for example, \cite{kolm2003acha}) that
\[ \int_{-\pi}^{\pi}\ln\left(\sin^2(s/2)\right)\cos(ks)ds=-\frac{2\pi}{|k|},\quad
\frac{1}{2\pi}\int_{-\pi}^{\pi}\cot\left(\frac{s}{2}\right)\sin(ks)ds=\operatorname{sgn}(k),
\]
we obtain
\be
a_k \approx \frac{1}{2}-\frac{r^2}{4\Delta t |k|^3},\quad
b_k \approx \frac{i}{2}\operatorname{sgn}(k)-\frac{ir^2}{4\Delta t k^2}\operatorname{sgn}(k).
\ee
Combining all the above, we see that for the Fourier mode $e^{iks}$ with $k$ large,
the following $2\times 2$ matrix determines its spectral behavior:
\be
\begin{bmatrix}\frac{1}{2}& \frac{i}{2}\operatorname{sgn}(k)[1-\frac{r^2}{2\Delta t k^2}]\\
  -\frac{i}{2}\operatorname{sgn}(k) & \frac{1}{2}\end{bmatrix}.
\ee
The above matrix has roughly equal eigenvalues and singular values with
 $\lambda_{k1}\approx\sigma_{k1}\approx 1-\frac{r^2}{8\Delta t k^2}$
and $\lambda_{k2}\approx\sigma_{k2}\approx \frac{r^2}{8\Delta t k^2}$ for $k$ large.
In summary, the integral equation system \eqref{bie} has eigenvalues (and roughly equal singular values)
$0$, $\lambda_0$, $\frac{r^2}{8\Delta t k^2}$, $1-\frac{r^2}{8\Delta t k^2}$
for $k$ large. From this, in the complement of the one-dimensional nullspace, 
the condition number of the linear system can be seen to be of the order 
$O(\Delta t/h^2)$ with $h =2\pi/N$ the spatial 
discretization size, assuming a uniform grid (so that $k \approx N \approx 1/h$).
In short, the condition number is
$O(N)$ for $\Delta t =O(h)$ and $O(1)$ for 
$\Delta t = O(h^2)$. For a fixed time step  $\Delta t$ independent of $N$, 
the condition number is
$O(N^2)$ (see \cref{condfig}). These estimates are
essentially the same as the conditioning of an implicit finite difference 
approximation applied to the heat equation with the same $\Delta t$ and $h$.

The preceding analysis can be extended, in part,
to the case of an arbitrary smooth curve.

\begin{lemma}\label{lem1}
The nullspace of the system of integral equations \eqref{bie} contains functions 
of the form $[\rho_0(\bx)f(t)\, \, \, 0]^T$ where $\rho_0(\bx)$ spans
the one-dimensional nullspace
of the operator $\frac{1}{2}I+S_{\rm L\nu}$ and 
$f(t)$ an arbitrary smooth function on $[0,T]$.
\end{lemma}

\noindent {\em Proof:}
  It is well known that the operator $\frac{1}{2}I+S_{\rm L\nu}$
  has a one-dimensional nullspace (see, for example, \cite{kress2014}). Let us
  denote a corresponding null vector by
  $\rho_0$. 
  It is easy to see that the function $v=\cs_{\rm L}[\rho_0]$ solves the interior
  Neumann problem for the Laplace equation with zero boundary data. 
  From well-known properties of harmonic functions, this implies
  that $v$ must be constant in $D$, so that its tangential derivative
  must be zero on the boundary. Thus, $S_{\rm L\tau}[\rho_0]=0$, completing the proof.
\qed
\begin{remark}
Numerical experiments indicate that the 
only null vectors of \eqref{bie} are the functions
identified in \eqref{lem1}.
We conjecture that the coupled system of integral equations 
is exactly rank one deficient in any simply connected domain.
\end{remark}

\section{Numerical results for the coupled integral equation system} \label{numsec}
Let us first consider the behavior of the
implicit, second-order accurate one-step marching
scheme described above. We will refer to solving the resulting system
of the form \eqref{syscircle} as the fully implicit combined source
integral equation (FI-CSIE).
As discussed in \cref{discret-sec},
we use a 16th order accurate
spatial quadrature rule \cite{alpert1999sisc} so that the 
spatial error is negligible and we accelerate the computation of 
the history part using the Fourier spectral method of 
\cite{greengard1990cpam,jiang2012sisc}.

\begin{figure}[ht]
\centering
\includegraphics[width=.90\textwidth]{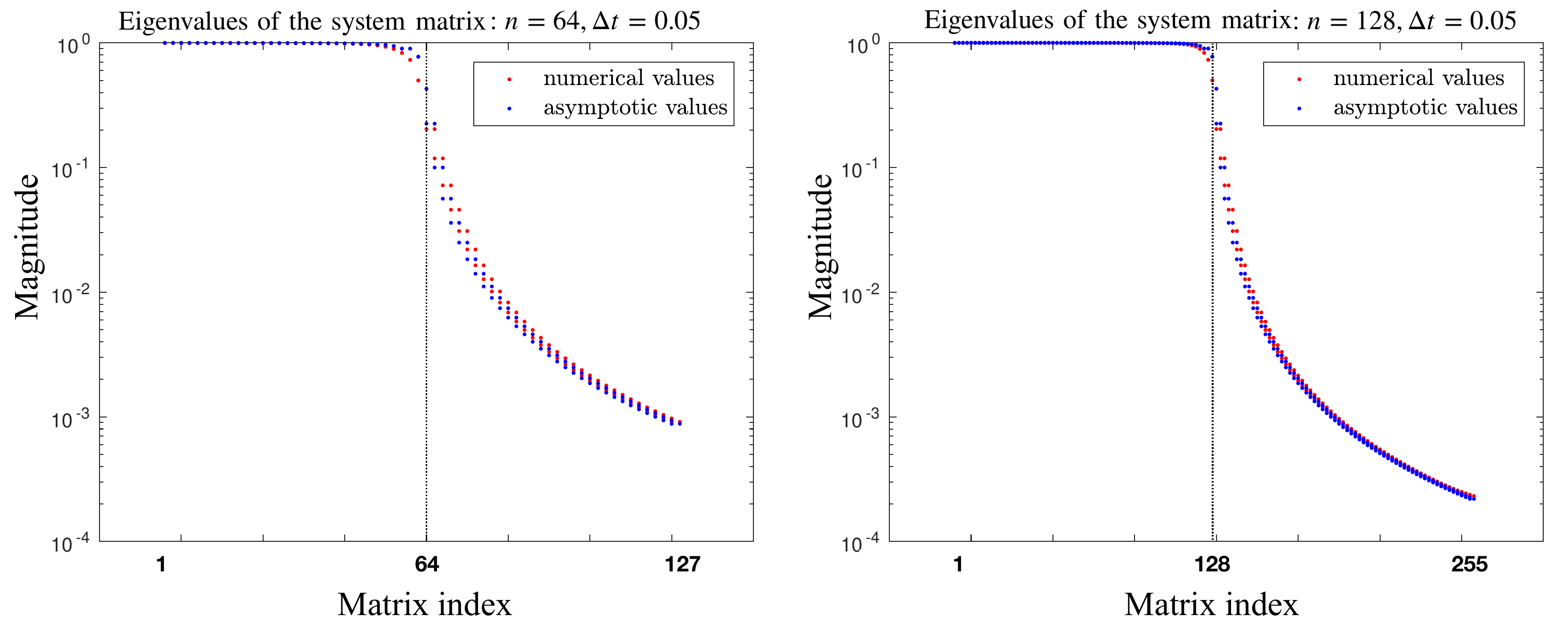}
\caption{Magnitude of the eigenvalues of the system matrix (of size $2n$)
  for a circle of radius $r=0.6$, discretized with $n$ points.
  The red dots are numerical values and the blue dots are the asymptotic values
  from \cref{specsec}. The $x$-axis corresponds to the eigenvalue index, 
  plotted in decreasing 
  order. The eigenvalues from the asymptotic analysis are ordered in the corresponding
  fashion.
  We omit the exact zero eigenvalue from the Fourier analysis, which is manifested by a 
  single eigenvalue of order $10^{-15}$ in the matrix analysis.}
\label{fig1}
\end{figure}

\begin{figure}[ht]
\centering
  \includegraphics[height=41mm]{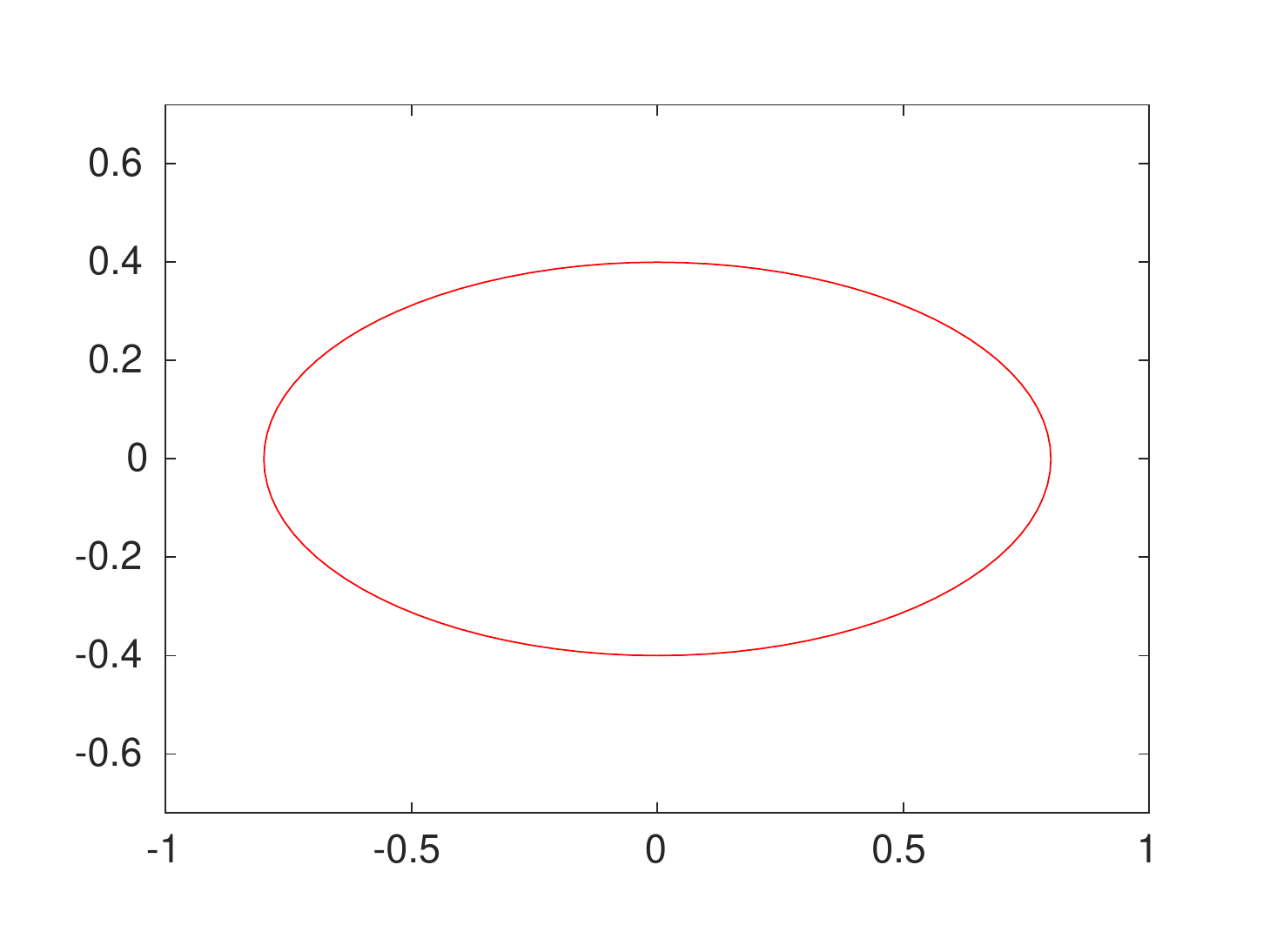}
\hspace{5mm}
  \includegraphics[height=40mm]{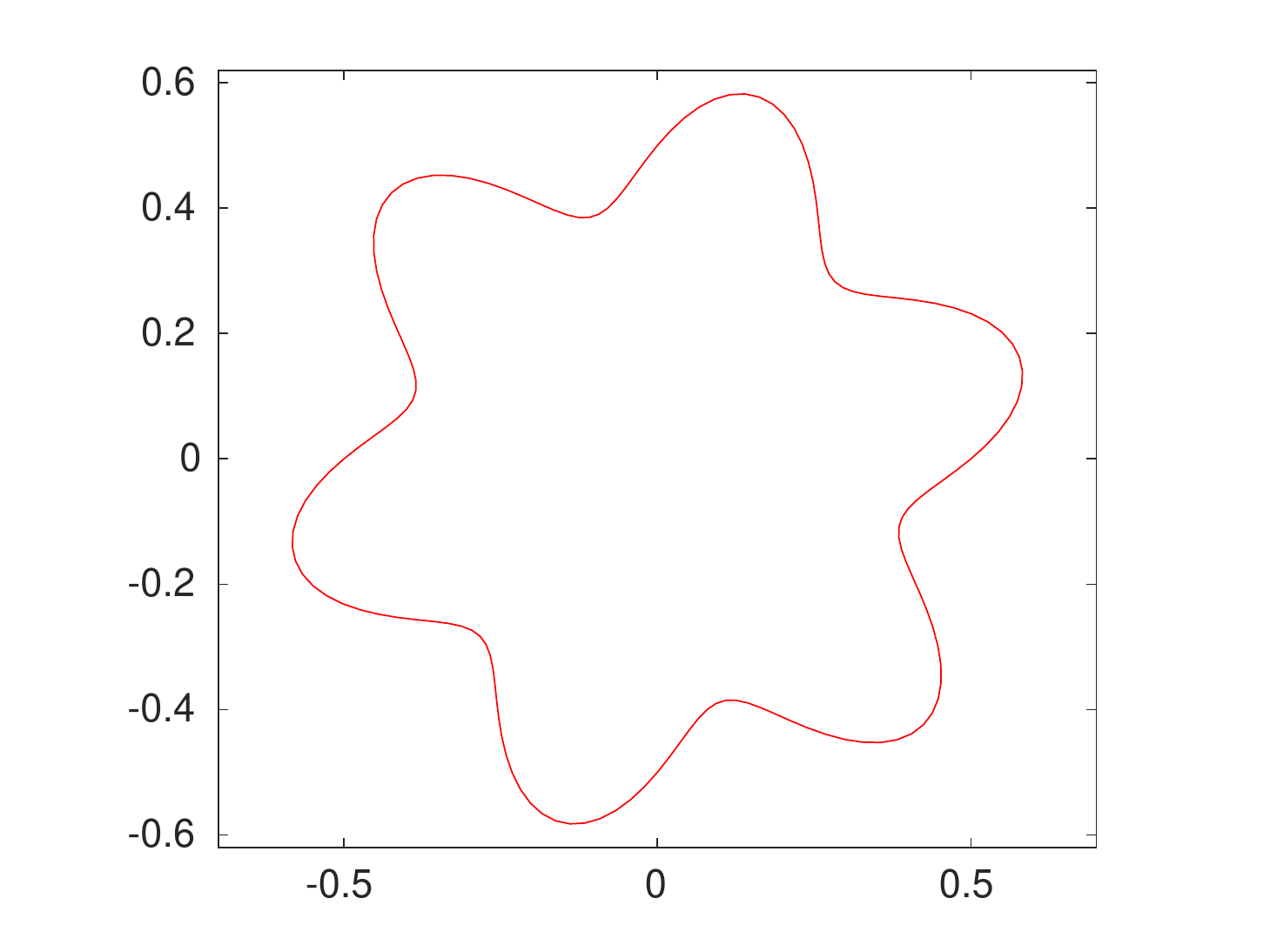}
\caption{Boundary curves for \cref{tab1,tab2}, respectively.}
\label{fig3}
\end{figure}

In \cref{fig1}, we plot the eigenvalues
of the system matrix with $n=64$ and $n=128$, respectively, when $\Gamma$ is
a circle of radius $r=0.6$. 
Note that the asymptotic analysis is in close agreement with direct discretization.
It is worth noting that a physical constraint on the boundary data ${\bf g}$ or 
$\tg$ is that the normal component integrates to zero on $\Gamma$. 
Such data has no projection onto the nullvector of the system matrix so that
an iterative method such as GMRES can be applied 
without any modification.
With a stopping criterion for the residual set to $10^{-12}$,
 \cref{tab1,tab2} show
the performance of GMRES and the obtained error 
when the boundary is either an ellipse with aspect ratio $2:1$ or a smooth hexagram,
respectively (\cref{fig3}). For all the tables presented in this paper, we take as the exact solution
the divergence-free velocity field
\be\label{extsol}
\begin{aligned} {\bf u}(\bx,t) &= 
 \sum_{j=1}^{10}\sum_{k=0}^{[t/(2h)]} \frac{\left(x_2-x_{2j},x_{1j}-x_1\right)}{|\bx-\bx_j|^2}
   \left(e^{-\frac{|\bx-\bx_j|^2}{4(t-(2k+1)h)}}-e^{-\frac{|\bx-\bx_j|^2}{4(t-(2k+2)h)}}\right)\\
&\quad   + t\cos(313\pi t)(x_1,-x_2)
   + \frac{t^2}{4}\cos(233\pi t) e^{x_1} ( \cos x_2, -\sin x_2)\\
&\quad 
   +
2t\sin(299\pi t) e^{x_2} ( \cos x_1, \sin x_1),
\end{aligned}
\ee
where $h=0.1$, ${\bf x} = (x_1,x_2)$ and the $\{ {\bf x}_j \}$ are chosen to be 
equispaced on the unit circle, which encloses both domains of interest. 
$96$ points are used to discretize the ellipse and $160$ points
are used to discretize the hexagram. In both tables, the first column lists
the total number of time steps $N$ needed to reach $t = 1$; 
the second column lists the time step;
the third column lists the average number of GMRES iterations required to solve the system
to the desired tolerance;
the fourth column lists the relative $l2$ error at 20 
random points in the computational domain; the last column lists
the ratio of the errors for each doubling of $N$.
Note that the data are consistent with second order accuracy 
in time. 
Note also that
the number of iterations is approximately equal to the number of points on
the boundary, as expected for a large time step with $\Delta t \approx 1/N$.

\begin{table}[ht]
  \begin{center}
{\bf \caption{Numerical results for the ellipse using the second order accurate
fully implicit marching scheme. The matrix size is $192\times 192$.
\label{tab1}
}}
\begin{tabular*}{1.0\textwidth}{@{\extracolsep{\fill}}  l  l  c l r  }
    \hline
    $N$& $\dt$ & $N_{\text{its}}$ & Error & Ratio \\ \hline 
10& $1/10$ & 93 & $2.0 \cdot 10^{-2}$ & \\
20& $1/20$ & 93 & $1.2 \cdot 10^{-2}$ & 1.6\\
40& $1/40$ & 92 & $3.8 \cdot 10^{-3}$ & 3.8\\
80& $1/80$ & 90 & $7.5 \cdot 10^{-4}$ & 4.2\\
160& $1/160$ & 86 & $1.8 \cdot 10^{-4}$ & 4.1\\
\hline
\end{tabular*}
\end{center}
\end{table}
\begin{table}[ht]
  \begin{center}
{\bf \caption{Numerical results for the hexagram using the second order accurate
fully implicit marching scheme. The matrix size is $320\times 320$.
\label{tab2}
}}
\begin{tabular*}{1.0\textwidth}{@{\extracolsep{\fill}}  l  l  c l r  }
    \hline
    $N$& $\dt$ & $N_{\text{its}}$ & Error & Ratio \\ \hline 
10& $1/10$ & 172 & $1.8 \cdot 10^{-2}$ & \\
20& $1/20$ & 178 & $6.9 \cdot 10^{-3}$ & 2.6\\
40& $1/40$ & 183 & $1.5 \cdot 10^{-3}$ & 4.5\\
80& $1/80$ & 185 & $2.8 \cdot 10^{-4}$ & 5.5\\
160& $1/160$ & 183 & $6.1 \cdot 10^{-5}$ & 4.6\\
\hline
\end{tabular*}
\end{center}
\end{table}
\section{The predictor corrector combined source integral equation (PC-CSIE)} 
\label{altsec}

The previous section shows that the fully implicit combined source integral
equation (FI-CSIE)
yields a somewhat ill-conditioned system of equations for large time steps. 
We now investigate a method for solving the unsteady Stokes equations using a
rule of predictor-corrector type.  

\begin{definition} \label{ie1}
The $k$th order predictor-corrector scheme for the CSIE,
denoted by PC-CSIE($k$), is given by 
\[ 
\ba
&{\rm Step\ 1:}\ {\rm Set}\ \mu_j = P^{E}_k[\bmu_{j-1}](j \Delta t),   \\
&{\rm Step\ 2:}\ {\rm Solve}\ \left( \frac{1}{2}I+ \cs_{\rm L\nu}  \right)
\rho_j = 
\nu \cdot \tg + \cs^{\rm far}_{\rm H\tau}[\bmu_{j-1}] + \cs^{\rm loc}_{\rm H\tau}[\bmu_{j}] , \\
&{\rm Step\ 3:}\ {\rm Solve}\ \left( \frac{1}{2}I+A^{j,0}_{\rm H\nu}  \right) \mu_j = 
\tau \cdot \tg- \cs^{\rm far}_{\rm H\nu}[\bmu_{j-1}] 
 - F_{\rm H\nu}^{j,k}[\bmu_{j-1}]
- \cs_{\rm L\tau}[\rho_{j}],
\ea
\]
using the notation of \cref{discret-sec}.
\end{definition}

That is, we first
extrapolate $\mu_j$ from previous time steps, then solve for $\rho_j$,
and finally solve for $\mu_j$ given the newly computed $\rho_j$.
Note that the integral equation in step 2 of PC-CSIE($k$)
is a Fredholm equation of the second kind.  While it has a one-dimensional
nullspace, it is simply the classical equation for the harmonic interior Neumann
problem obtained when representing the solution as a single layer potential. 
Since the right-hand side is compatible,
the resulting linear system is easily solved using GMRES with $O(1)$ iterations.
The integral equation in step 3 is a Volterra equation of the second
kind. It is always invertible and well-conditioned (at least
on smooth curves).

The obvious difference between PC-CSIE ($k$) 
and the fully implicit version comes from the extrapolation step. 
While the order
of accuracy can be arbitrarily high, and the individual integral equations are
well-conditioned, the stability of the resulting scheme remains to be studied.

Preliminary experiments suggest that the schemes PC-CSIE(2) and 
PC-CSIE(3) are stable even for moderately large time steps, while 
PC-CSIE(4) is not, but a thorough analysis
remains to be undertaken (see \cref{sec:conc}).
Using the same geometries 
as above (\cref{fig3}), we obtain the results in 
\Cref{tab3} for the ellipse (top) and the hexagram (bottom).
Note that many fewer GMRES iterations are required 
for each step. The convergence rates, estimated by the error ratios
in the last column are somewhat erratic, but generally better than
the theoretical estimate $2^{k}$.
\begin{table}[ht]
\begin{center}
{\bf \caption{Numerical results for the ellipse and hexagram using the 
PC-CSIE(2) method. 
The matrix size is $96\times 96$ for the ellipse
and $256\times 256$ for the hexagram at each stage.
\label{tab3}
}}
\begin{tabular*}{1.0\textwidth}{@{\extracolsep{\fill}}  l  l  c l r  }
\hline\hline
Ellipse\\
\hline
$N$& $\dt$ & $N_{\text{its}}$ & Error & Ratio \\ \hline 
80& $1/40$ & 7 & $4.2 \cdot 10^{-2}$ & \\
160& $1/80$ & 7 & $1.7 \cdot 10^{-3}$ & 24\\
320& $1/160$ & 7 & $2.9 \cdot 10^{-4}$ & 6.0\\
640& $1/320$ & 7 & $4.8 \cdot 10^{-5}$ & 6.0\\
1280& $1/640$ & 7 & $9.0 \cdot 10^{-6}$ & 5.3\\
\hline\hline
Hexagram\\
\hline
$N$& $\dt$ & $N_{\text{its}}$ & Error & Ratio \\ \hline 
80& $1/40$ & 12 & $9.6 \cdot 10^{-3}$ & \\
160& $1/80$ & 12 & $3.6 \cdot 10^{-3}$ & 2.7\\
320& $1/160$ & 11 & $5.8 \cdot 10^{-4}$ & 6.1\\
640& $1/320$ & 11 & $9.9 \cdot 10^{-5}$ & 5.8\\
1280& $1/640$ & 11 & $1.7 \cdot 10^{-5}$ & 5.7\\
\hline\hline
\end{tabular*}
\end{center}
\end{table}

\section{High order schemes using spectral deferred corrections} \label{sdcsec}

Instead of seeking to develop stable higher order predictor-corrector type schemes,
we now show that more rapid convergence is easily achieved by combining 
PC-CSIE(2) with
spectral deferred correction (SDC)
\cite{blmsdc,dutt2000bit,hagstromsdc,minionsdc}.

A very brief introduction to deferred corrections follows:
suppose that we seek the solution $v(t)$ of some 
time-dependent problem starting at $t=0$, 
and that an approximate solution can be computed for $k$ 
steps on $[0,\Delta t]$ 
using a low order accurate method, with an error 
of the order $O(\Delta t^m)$ for some $m<k$. We denote the discrete solution at those
$k$ points by ${\bf v}^{[0]}_k$. One can then interpolate the low order solution
by a polynomial in $t$ of order $k$, namely $P_k^I[{\bf v}^{[0]}_k]$, 
defined in \cref{discret-sec}.
This allows us to defines a continuous error function
\[ \delta^{[0]}(t) = v(t) - P_k^I[{\bf v}^{[0]}_k](t), \] 
which can be substituted into the governing 
equation for $v(t)$ and solved for $\delta^{[0]}(t)$, using the same low order scheme.
This generates the discrete solution vector ${\bdelta}^{[0]}_k$.
A corrected approximation is then defined by
\[ {\bf v}^{[1]}_k = {\bf v}^{[0]}_k + {\bdelta}^{[0]}_k. \] 
It is straightforward to show that the error in  
${\bf v}^{[1]}_k$ is of the order
$O(h^{2m})$, so long as $2m < k$ and all computations involving the known function
$P_k^I[{\bf v}^{[0]}_k](t)$ are carried out with $k$th order accuracy. For further
details, see the references above and \cite{dcbook}. The correction procedure
is easily iterated until $k$th order accuracy is achieved. The process can then be
repeated on the next time interval $[\Delta t, 2\Delta t]$, etc.
(The phrase {\em spectral deferred correction} is typically used when the underlying 
problem has been formulated as an integral equation and the
$k$ stages are chosen at nodes corresponding to some high order
spectral discretization, typically of Gauss or Gauss-Radau type.)

In the present context,
let us assume that we have divided
$[0,T]$ into
$N$ equal subintervals $[t_{i-1},t_{i}]$ with $t_i=i\dt$, $\dt=T/N$ for
$i=1,\ldots,N$. We restrict our attention to the $i$th 
such interval $[t_{i-1},t_i]$ which we will denote by $[\alpha,\beta]$ when the
context is clear.
Given a positive integer $k$, we will denote by 
$\alpha < \tau_1, \ldots, \tau_k=\beta$ the $k$ Gauss-Radau nodes shifted and scaled
to the interval $[\alpha,\beta]$ (see, for example, \cite{gautschi2000}),
with $\tau_0 = \alpha$.
Let us denote the values of the densities $\rho$ and $\mu$
at these nodes 
by ${\brho}_k=(\rho_0,\rho_1,\ldots,\rho_k)^T$ and
${\bmu}_k=(\mu_0,\mu_1,\ldots,\mu_k)^T$, respectively. 
Since we are discretizing in time only, recall that $\rho_i,\mu_i$ are functions
of the spatial variable $\by\in \Gamma(t)$.

Following the principle outlined above,
the first step of SDC for the mixed potential formulation
is to use some low order scheme to obtain ${\brho}_k^{[0]}$ and 
${\bmu}_k^{[0]}$.
We then use these two vectors to obtain interpolating polynomials
of degree $k-1$ in time, namely
$P_k^I[{\bmu}^{[0]}_k]$ and 
$P_k^I[{\brho}^{[0]}_k]$ and define 
\[
\delta_{\mu}^{[0]}(t) = \mu(t) - P_k^I[{\bmu}^{[0]}_k],\quad
\delta_{\rho}^{[0]}(t) = \rho(t) - P_k^I[{\brho}^{[0]}_k].
\]
Inserting this representation into \cref{bie}, we obtain

\be\label{biecorr}
\ba
\frac{1}{2}\delta^{[0]}_\rho(\bx,t)+
\cs_{\rm L\nu}[\delta^{[0]}_\rho](\bx,t)+\cs_{\rm H\tau}[\delta^{[0]}_\mu](\bx,t)&=
R_1(\bx,t), \\
-\cs_{\rm L\tau}[\delta^{[0]}_\rho](\bx,t)+\frac{1}{2} \delta^{[0]}_\mu(\bx,t)+
\cs_{\rm H\nu}[\delta^{[0]}_\mu](\bx,t)&=
R_2(\bx,t),
\ea
\ee
where the residuals $R_1$ and $R_2$ are given by
{\small
\[
\ba
R_1(\bx,t) &= \nu \cdot \tg(\bx,t) -
\frac{1}{2} P_k^I[{\brho}^{[0]}_k](\bx,t)-
\cs_{\rm L\nu}[P_k^I[{\brho}^{[0]}_k]](\bx,t)-
\cs_{\rm H\tau}[P_k^I[{\bmu}^{[0]}_k]](\bx,t),\\
R_2(\bx,t) &= 
-\tau \cdot \tg(\bx,t) +
\cs_{\rm L\tau}[P_k^I[{\brho}^{[0]}_k]](\bx,t)-
\frac{1}{2} P_k^I[{\bmu}^{[0]}_k](\bx,t)-
\cs_{\rm H\nu}[P_k^I[{\bmu}^{[0]}_k]](\bx,t).
\ea
\]
}

Note that this is exactly the same equation as \eqref{bie}, but with a different
right-hand side. Since, as noted above, SDC requires that all residuals 
be computed with high order accuracy, we have provided some of the integrals 
needed at the intermediate stages in the Appendix.

After solving \cref{biecorr} at the same $k$ stages, yielding 
$[{\bdelta}^{[0]}_\rho]_k$, $[{\bdelta}^{[0]}_\mu]_k$,
we let
\be
{\brho}^{[1]}_k = 
{\brho}^{[0]}_k +
[{\bdelta}^{[0]}_\rho]_k, \quad 
{\bmu}^{[1]}_k = 
{\bmu}^{[0]}_k +
[{\bdelta}^{[0]}_\mu]_k.
\ee
This procedure may be repeated until the desired order of accuracy is achieved.

We have implemented SDC using the second order predictor corrector
scheme PC-CSIE(2) described in the previous section. 
For simplicity, we provide numerical results
in \cref{tab4} for the case of a circle. In this table, 
$\text{SDC}_k^j$ denotes the scheme with $k$ Gauss-Radau nodes on each
subinterval and $j$ iterations of deferred correction. 
In particular, 
$\text{SDC}_k^0$ is simply the uncorrected solution obtained with
$\text{PC-CSIE(2)}$.
$N$ is
the number of subintervals, $\dt$ is the time step size for each subinterval,
$E$ is the relative $l^2$ error for the method indicated in the subscript, and 
${N}_{\text{it}}$ is the average number of iterations for GMRES to 
reach the requested tolerance $10^{-12}$. The
total number of time steps is $Nk$ and the expected error reduction should be
$2^{2(j+1)}$ for $\text{SDC}_k^{j}$ until $2(j+1) > k$, 
since we are driving the deferred correction process with a second order accurate
scheme.
\begin{table}[ht]
  \begin{center}
{\bf \caption{Numerical results for the circle of radius 0.5, 
using spectral deferred
correction. The number of points in the spatial discretization is $200$.
\label{tab4}
}}
{\small 
\begin{tabular*}{1.0\textwidth}{@{\extracolsep{\fill}}  |l|c|c|l|l|l|l|l|}
  \hline
 $Nk$ & $\dt/k$ & ${N}_{\text{its}}$ & $E_{SDC^0_5}$ & $E_{SDC^1_5}$ & 
$E_{SDC^2_5}$ & $E_{SDC^3_5}$ & $E_{SDC^4_5}$ 
    \\ \hline 
$40$ & $1/20$ & $2.7$ & $1.5\cdot 10^{-2}$ & $2.5\cdot 10^{-2}$ &
$2.5\cdot 10^{-2}$ & $2.5\cdot 10^{-2}$ & $2.5\cdot 10^{-2}$  \\
$80$ & $1/40$ & $2.6$ & $3.0\cdot 10^{-3}$ & $3.5\cdot 10^{-5}$ &
$1.8\cdot 10^{-4}$ & $1.8\cdot 10^{-4}$ & $1.8\cdot 10^{-4}$  \\
$160$ & $1/80$ & $2.7$ & $2.7\cdot 10^{-3}$ & $3.8\cdot 10^{-5}$ &
$3.7\cdot 10^{-5}$ & $3.7\cdot 10^{-5}$ & $3.7\cdot 10^{-5}$  \\
$320$ & $1/160$ & $2.6$ & $6.9\cdot 10^{-4}$ & $1.5\cdot 10^{-7}$ &
$1.3\cdot 10^{-7}$ & $9.8\cdot 10^{-8}$ & $4.5\cdot 10^{-8}$  \\
$640$ & $1/320$ & $2.6$ & $1.6\cdot 10^{-4}$ & $1.1\cdot 10^{-8}$ &
$5.3\cdot 10^{-9}$ & $4.5\cdot 10^{-9}$ & $4.9\cdot 10^{-9}$  \\
\hline
\end{tabular*}
}
\end{center}
\end{table}
While the behavior of the SDC schemes as a function of the number of 
correction sweeps is somewhat erratic, it is more or less consistent 
with the asymptotic estimates.
That is, $SDC^0_5$ converges approximately like a second order scheme,
while $SDC^1_5$ converges at a much higher rate.
Further sweeps of deferred correction don't increase the convergence
rate significantly, since the degree of polynomial approximation 
is only four, limiting the order of accuracy, as discussed above.
Note, however, that these further sweeps have no impact on
stability.
\begin{remark}
A nice feature of the mixed potential representation is the 
complete separation of the {\underline {\em instantaneous}} 
pressure source from the vortex source.
As a result, even though the
exact solution defined in \cref{extsol} has a highly
oscillatory pressure field, the numerical results
in \cref{tab1,tab2,tab3,tab4} show
that high accuracy is achieved even for large time steps that 
under-resolve the oscillatory behavior of the pressure.
\end{remark}
\section{Conclusions and future work} \label{sec:conc}

We have developed a new integral representation for the unsteady
Stokes equations which makes use of simple harmonic and heat 
potentials, leading to the {\em combined source integral equation} 
(CSIE). Unlike schemes based on the unsteady Stokeslet
\cite{jiang2012sisc}, this permits 
the direct application of well-developed fast algorithms (see, for
example,
\cite{greengard2000acha,greengard1990cpam,lubichheat1,tauschheat1,wang2018,
fullheatsolver}
and
\cite{fmm2,ethridge2001sisc,GLee,fmm,langston2011,biros2015cicp}).

While the fully coupled CSIE is not of Fredholm
type, we have shown that each individual equation is well-conditioned, and
found that a second-order predictor-corrector type method is effective even
for large time steps.  Moreover, one can achieve high order
accuracy through the use of spectral deferred correction.
Since our primary goal in the present paper is the development of the 
mathematical
representation itself, a more thorough investigation of various
predictor-corrector, Runge-Kutta, and implicit-explicit type marching schemes
will be carried out at a later date.

It is worth noting that the mixed potential representation and existing 
gauge methods (see \cref{gauge}) bear some resemblance.
The principle differences are (1) that we are working in an integral 
equation-based
framework, (2) that we apply the Helmholtz decomposition to the 
inhomogeneous data rather than to the auxiliary unknown vector
field $\bm$, and (3) that we have as unknowns only the vortex source and 
pressure source, which are restricted to the boundary and for which imposing 
velocity boundary conditions is straightforward.

A number of open questions remain, including the 
completeness of the representation for multiply-connected domain,
a detailed characterization of the nullspace of the coupled system,
and the extension of the mixed potential formulation to other 
boundary/interface conditions.  While our representation is formally 
valid in either fixed or moving geometries, 
we have not yet investigated its performance in the nonstationary case.
To solve the equations with a forcing term 
(or the full Navier-Stokes equations),
we also need to couple the solver described here with volume-integral based 
codes for the Helmholtz decomposition, as discussed in section \ref{helmsec}.
We are presently investigating all of these topics and will report on our
progress at a later date.

\section*{Appendix} \label{app1}

For the SDC method of \cref{sdcsec}, using
product integration in time as in \cite{jiang2012sisc}, 
the residual is required at intermediate stage $\tau_i \in [\alpha,\beta]$.
This involves integrals beyond those given by eq.~(4.9)
in \cite{jiang2012sisc}. We provide those integrals here, which are sufficient
to obtain fourth order accuracy.

\be
\begin{aligned}
  &\int_{\alpha}^{\alpha+\tau_i} e^{-r^2/4(\alpha+\tau_i-\tau)}
  \frac{(\beta-\tau)^j}{(\alpha+\tau_i-\tau)^2}d\tau \\
&=\left\{\begin{aligned}&\frac{4}{r^2}\ec, & j=0,\\
& E_1(c)+\frac{b}{a}\ec, & j=1,\\
&\left(\tau_i+\frac{b^2}{a}\right)\ec-E_1(c)(a-2b), & j=2,\\
&-\frac{\ec}{2}\left(\tau_ia-x_i^2-6\tau_ib-2b^3/a \right)+E_1(c)(a^2-6ab+6b^2), & j=3,\\
&\frac{\ec}{6}\left(2\tau_i^3-\tau_i^2(a-12b)+\tau_i(a-6b)^2+6b^4/a)\right)& \\
&-\frac{E_1(c)}{6}(a^3-12a^2b+36ab^2-24b^3), & j=4,\\
\end{aligned}
\right.
\end{aligned}
\ee
where $r=\|\bx-\by\|$, $c=\frac{r^2}{4\tau_i}$, $a=\frac{r^2}{4}$, $b=\dt-\tau_i$.

\section*{Acknowledgments}

We would like to thank Alex Barnett, Charlie Epstein, Jingfang Huang, Manas Rachh,
Shravan Veerapaneni and Jun Wang for many useful conversations.

\bibliographystyle{siam}
\bibliography{all2}

\end{document}